\documentclass[reqno]{amsart}

\usepackage{easywncy}

\newcommand{\shc}{\widetilde{\mathcyr {sh}}}
\newcommand{\sh}{{\mathcyr {sh}}}
\newcommand{\s}{\star}
\newcommand{\Z}{{\mathfrak Z}}
\newcommand{\h}{{\mathfrak H}}
\newcommand{\sym}{{\mathfrak S}}

\newcommand{\Q}{{\mathbb Q}}
\newcommand{\reg}{{\rm reg}_{\sh}}
\theoremstyle{plain}
\newtheorem{thm}{Theorem}[section]
\newtheorem{conj}[thm]{Conjecture}

\newtheorem{prop}[thm]{Proposition}

\theoremstyle{definition}

\newtheorem{exmp}[thm]{Example}

\begin{document}

\markboth{K. Imatomi, T. Tanaka, K. Tasaka and N. Wakabayashi}
{On some combinations of MZSV's}

%
%

\title{On some combinations of multiple zeta-star values
}

\author{Kohtaro Imatomi
}

\address{Graduate School of Mathematics, Kyushu University \\
744 Motooka Fukuoka-city, Fukuoka, 819-0395 Japan,
}
\email{k-imatomi@math.kyushu-u.ac.jp
} 

\author{Tatsushi Tanaka
}

\address{Faculty of Mathematics, Kyushu University \\
744 Motooka Fukuoka-city, Fukuoka, 819-0395 Japan,
}
\email{t.tanaka@math.kyushu-u.ac.jp
}

\author{Koji Tasaka
}

\address{Graduate School of Mathematics, Kyushu University \\
744 Motooka Fukuoka-city, Fukuoka, 819-0395 Japan,
}
\email{ma209026@math.kyushu-u.ac.jp
}

\author{Noriko Wakabayashi
}

\address{Faculty of Engineering, Kyushu Sangyo University \\
3-1 Matsukadai 2-chome, Higashi-ku, Fukuoka, 813-8503 Japan,
}
\email{noriko@ip.kyusan-u.ac.jp
}

\maketitle


\begin{abstract}
We prove that the sum of multiple zeta-star values over all indices inserted two $2$'s into the string $3,1,\ldots ,3,1$ is evaluated to a rational multiple of powers of $\pi^2$. We also establish certain conjectures on evaluations of multiple zeta-star values observed by numerical experiments. 
\end{abstract}



\section{Introduction/Main theorem}
We discuss some explicit evaluations of multiple zeta values (MZV's for short) and multiple zeta-star values (MZSV's for short), 
which are defined, for positive integers $k_1, k_2, \ldots ,k_n $ with $k_1 \geq 2$, by the convergent series
\[ \zeta (k_1 ,k_2, \ldots , k_n) = \sum_{m_1>m_2> \cdots >m_n>0} \frac{1}{m_1^{k_1} m_2^{k_2} \cdots m_n^{k_n}} \]
and
\[ \zeta^{\s} (k_1 ,k_2, \ldots , k_n) = \sum_{m_1\geq m_2\geq  \cdots \geq m_n\geq 1} \frac{1}{m_1^{k_1} m_2^{k_2} \cdots m_n^{k_n}}, \]
respectively. The number $k_1+k_2+\cdots +k_n$ is called weight
. When $n=1$, MZV's and MZSV's coincide and are known as the special value of the Riemann zeta function at positive integers. 

Special values of the Riemann zeta function have been investigated since Euler. Several results are found for example in \cite{Euler,apery,rivoal,zdl}. For MZV's, the following evaluation was proved by many
authors \cite[$\ldots$]{H,don,BBB,yamasaki}: for $m>0$, we have
\[ \zeta ( \underbrace{ 2,\ldots ,2}_{m} )=\frac{ \pi^{2m}}{(2m+1)!}. \]
Zagier conjectured the formula
\[ \zeta ( \underbrace{3,1, \ldots ,3,1}_{2n} ) = \frac{2\pi^{4n}}{(4n+2)!} \]
for $n>0$ in \cite{don}, which was solved in 
\cite{BBBL1} by using certain property of the iterated integral shuffle product rule. Kontsevich and Zagier
gave another proof of the formula in connection with the Gauss hypergeometric function (\cite{KZ}). 
In \cite{BB,muneta2}, the following more general identity is proved: for any non-negative integers $n$ and $m$ with $n+m>0$, we have
\begin{align} \label{bb} \sum_{\begin{subarray}{c} j_0+j_1+\cdots +j_{2n}=m \\ j_0,j_1,\ldots ,j_{2n} \geq 0 \end{subarray}} \zeta ( \{ 2 \}^{j_0} , 3, \{ 2 \}^{j_1} ,1,\ldots , 3, \{ 2 \}^{j_{2n-1}} ,1 ,\{ 2 \}^{j_{2n}} ) \nonumber \\
= \binom{ m+2n}{m } \frac{\pi^{2m+4n}}{(2n+1)(2m+4n+1)!}, 
\end{align}
where $\{ 2 \}^j$ stands for the $j$-tuple of $2$. More precise conjecture is introduced in \cite[Conjecture $1$]{BBBL1}. 

While the special evaluations are well studied in the MZV's case, there are less observations on special evaluations of MZSV's. For $m>0$, the property
\[ \zeta^{\s} ( \underbrace{ 2,\ldots ,2}_{m} )\in\Q\cdot\pi^{2m} \]
is proved for example in \cite{H,zlob,ako,yamasaki}. In \cite{muneta1}, formulas\[ \zeta^{\s} ( \underbrace{3,1, \ldots ,3,1}_{2n})\in\Q\cdot\pi^{4n} \]
for $n>0$ and
\[ \sum_{\begin{subarray}{c} j_0+j_1+\cdots +j_{2n}=1 \\ j_0,j_1,\ldots ,j_{2n} \geq 0 \end{subarray}} \zeta^{\s} ( \{ 2 \}^{j_0} , 3, \{ 2 \}^{j_1} ,1,\ldots , 3, \{ 2 \}^{j_{2n-1}} ,1 ,\{ 2 \}^{j_{2n}} ) \in \Q \cdot \pi^{4n+2} \]
for $n\geq 0$ are proved. 

In the present paper, we prove the following result. 
\begin{thm} \label{main} For any non-negative integer $n$, we have
\begin{align*} \sum_{\begin{subarray}{c} j_0+j_1+\cdots +j_{2n}=2 \\ j_0,j_1,\ldots ,j_{2n} \geq 0 \end{subarray}} \zeta^{\s} ( \{ 2 \}^{j_0} , 3, \{ 2 \}^{j_1} ,1,\ldots , 3, \{ 2 \}^{j_{2n-1}} ,1 ,\{ 2 \}^{j_{2n}} ) \in \Q \cdot \pi^{4n+4}. 
\end{align*} 
\end{thm}
\noindent In section \ref{sec2}, we state a bit more general formula, which yields Theorem \ref{main}. In section \ref{sec3}, we prove our main theorem (Theorem \ref{main}). In section \ref{sec4} we give more precise conjectures on some special evaluations of MZSV's, which are observed by numerical experiments. We give some remarks on our conjectures in section \ref{sec5}.

\section{Key identities}\label{sec2}
To prove our main theorem, we describe some key identities here. It is suitable for stating our key identities to use the algebraic setup of MZV's and MZSV's introduced by Hoffman \cite{hoff} as follows. 

Let the symbol $\h$ be the noncommutative polynomial algebra in two indeterminates $x$ and $y$, and $\h^1 ,\h^0$ its subalgebras:
\[ \h := \Q \langle x,y \rangle \supset \h^1 := \Q + \h y \supset \h^0 := \Q +x \h y.\]
The degree of a word is called its weight. We put $z_l = x^{l-1} y$ for $l \geq 1$. The algebra $\h^1$ is the noncommutative polynomial algebra over $\Q$ freely generated by $\{ z_1, z_2, z_3, \ldots \}$. We define two $\Q$-linear maps $Z,\overline{Z}:\h^0 \to {\mathbb R}$ respectively by
\[ Z(z_{k_1} z_{k_2} \cdots z_{k_n} )=\zeta (k_1,k_2, \ldots ,k_n),\ Z(1)=1,   \]
and 
\[ \overline{Z}(z_{k_1} z_{k_2} \cdots z_{k_n} )=\zeta^{\s} (k_1,k_2, \ldots ,k_n),\ \overline{Z}(1)=1, \]
which are usually called the evaluation maps. The weight of a word is that of the corresponding MZV or MZSV. 

It is well known that any MZSV can be expressed as a $\Q$-linear combination of MZV's and vice versa. 
For example, we have
\begin{align} \zeta^{\s} (k_1, k_2) &= \sum_{n \geq m \geq 1} \frac{1}{n^{k_1} m^{k_2}} \nonumber \\
&= \sum_{n > m>0} \frac{1}{n^{k_1} m^{k_2}} + \sum_{n=m \geq 1} \frac{1}{n^{k_1} m^{k_2}} = \zeta (k_1 ,k_2) +\zeta (k_1+k_2). \label{e1} 
\end{align}
\noindent 
Let $\gamma$ be the automorphism on $\h$ characterized by
\[ \gamma (x)=x ,\ \gamma (y)=x+y. \]
We define the $\Q$-linear map $d: \h^1 \rightarrow \h^1 $ by
\[  d(wy)=\gamma (w)y \]
for any word $w \in \h$.
Then 
the linear transformation between MZV's and MZSV's is expressed as
\[ \overline{Z} = Z \circ d. \]

Let $\ast\colon\h^1 \times \h^1 \rightarrow \h^1$ be the $\Q$-bilinear map defined, for any words $w,w' \in \h^1$ and any positive integers $k_1,k_2$, by
\[ 1\ast w=w\ast 1=w \]
\noindent and the recursive rule
\begin{align} \label{har} z_{k_1} w\ast z_{k_2} w' =z_{k_1}(w\ast z_{k_2}w')+z_{k_2}(z_{k_1}w\ast w') +z_{k_1+k_2} (w\ast w'). 
\end{align}
It is known that the product $\ast$ is commutative and associative (\cite{hoff}). The product $\ast$ is called the harmonic product on $\h^1$. We find that $\h^0\ast\h^0 \subset \h^0$ and the map $Z$ is a homomorphism with respect to the harmonic product.

As is defined in \cite{muneta2}, we introduce another $\Q$-bilinear map $\shc:\h^1 \times \h^1 \rightarrow \h^1$ by
\[ 1 \ \shc \ w = w \ \shc \ 1 = w \]
\noindent and
\begin{align}
z_{k_1} w \ \shc \ z_{k_2} w^{\prime} = z_{k_1} (w \ \shc \ z_{k_2} w^{\prime} ) + z_{k_2} (z_{k_1} w \ \shc \ w^{\prime} ) \label{eqn:ast} 
\end{align}
for any words $w,w^{\prime} \in \h^1$ and any positive integers $k_1, k_2$. We see that the product $\shc$ is commutative and associative. However, we notice that each of the evaluation maps $Z$ and $\overline{Z}$ can not be a homomorphism with respect to the product $\shc$.

Under these notations, we have the following formulas.

\begin{thm} \label{mainthm} Let $a,b,c$ be positive integers. For any integer $n\geq 0$, we have
\begin{equation*}(\alpha_n)\ \begin{aligned}[t] &d(z_c^2 \ \shc \ (z_a z_b)^n) \\
&=2 \sum_{j+k=n} d(z_c \ \shc \ (z_a z_b)^j)\ast z_{(a+b)k+c} +\sum_{j+k=n} (z_c^2 \ \shc \ (z_a z_b)^j)\ast d(z_{a+b}^k) \\ 
&\ -4\sum_{i+j+k=n} d((z_a z_b)^i)\ast z_{(a+b)j+c} z_{(a+b)k+c} -\sum_{j+k=n} d((z_a z_b )^j)\ast z_{(a+b)k+2c}, 
\end{aligned}
\end{equation*}
\begin{equation*}(\beta_n)\ \begin{aligned}[t] &d(z_c^2 \ \shc \ z_b (z_a z_b)^n) \\
&=2 \sum_{j+k=n} d(z_c \ \shc \ z_b(z_a z_b)^j)\ast z_{(a+b)k+c}  +\sum_{j+k=n} (z_c^2 \ \shc \ z_b(z_a z_b)^j)\ast d(z_{a+b}^k) \\
&\ -4\sum_{i+j+k=n} d(z_b(z_a z_b)^i)\ast z_{(a+b)j+c} z_{(a+b)k+c} -\sum_{j+k=n} d(z_b(z_a z_b )^j)\ast z_{(a+b)k+2c}. 
\end{aligned}
\end{equation*}
\end{thm}
\noindent 
Theorem \ref{mainthm} is the core property to prove our main theorem (Theorem \ref{main}). The proofs of theorems are presented in the next section.

\section{Proofs}\label{sec3}
First we prove Theorem \ref{main} by assuming Theorem \ref{mainthm}, the proof of which is given next.
\begin{proof}[Proof of Theorem \ref{main}]
By putting $a=3,b=1$ and $c=2$ into $(\alpha_n)$ of Theorem \ref{mainthm}, we have
\begin{align}
\begin{split}
d(z_2^2 \ \shc \ (z_3 z_1)^n) &=2 \sum_{j+k=n} d(z_2 \ \shc \ (z_3 z_1)^j)\ast z_{4k+2}+\sum_{j+k=n} (z_2^2 \ \shc \ (z_3 z_1)^j)\ast d(z_4^k) \\
&\ -4\sum_{i+j+k=n} d((z_3 z_1)^i)\ast z_{4j+2} z_{4k+2}-\sum_{j+k=n} d((z_3 z_1 )^j)\ast z_{4k+2} . \label{a}
\end{split}
\end{align}
\noindent By the harmonic product rule \eqref{har}, the third term of the right-hand side of (\ref{a}) can be written as
\[  -2\sum_{i+j+k=n} d((z_3 z_1)^i)\ast (z_{4j+2} \ast z_{4k+2} - z_{4j+4k+4} ) \]
\noindent Evaluating (\ref{a}) via the map $Z$, we obtain
\begin{align}
\begin{split}
& \sum_{\begin{subarray}{c} j_0+j_1+\cdots +j_{2n}=2 \\ j_0,j_1,\ldots ,j_{2n} \geq 0 \end{subarray}} \zeta^{\s} ( \{ 2 \}^{j_0} , 3, \{ 2 \}^{j_1} ,1,\ldots , 3, \{ 2 \}^{j_{2n-1}} ,1 ,\{ 2 \}^{j_{2n}} )  \\
& =2 \sum_{i=0}^n \sum_{\begin{subarray}{c} j_0+j_1+\cdots +j_{2i}=1 \\ j_0,j_1,\ldots ,j_{2i} \geq 0 \end{subarray}} \zeta^{\s} ( \{ 2 \}^{j_0} , 3, \{ 2 \}^{j_1} ,1,\ldots , 3, \{ 2 \}^{j_{2i-1}} ,1 ,\{ 2 \}^{j_{2i}} )  \zeta (4n-4i+2) \\
&\ +\sum_{i=0}^n \sum_{\begin{subarray}{c} j_0+j_1+\cdots +j_{2i}=2 \\ j_0,j_1,\ldots ,j_{2i} \geq 0 \end{subarray}} \zeta ( \{ 2 \}^{j_0} , 3, \{ 2 \}^{j_1} ,1,\ldots , 3, \{ 2 \}^{j_{2i-1}} ,1 ,\{ 2 \}^{j_{2i}} )  \zeta^{\s} ( \{ 4 \}^{n-i} ) \\
&\ -2\sum_{i+j+k=n} \zeta^{\s} ( \{3,1\}^i ) \bigl\{ \zeta (4j+2) \zeta (4k+2) -\zeta (4j+4k+4) \bigr\} \\
&\ - \sum_{j+k=n} \zeta^{\s} ( \{ 3,1 \}^j ) \zeta (4k+2),
\end{split}\label{c}
\end{align}
\noindent
where $\{ 3,1 \}^l$ stands for the string $\underbrace{3,1,\ldots ,3,1}_{2l}$ and MZ(S)V of the empty index is regarded as $1$. We know $\zeta (2n) \in \Q \cdot \pi^{2n} $, $\zeta^{\s} (\{ 4 \}^n) \in \Q \cdot \pi^{4n}$, the formula (\ref{bb}) for $m=2$ and
\[ \sum_{\begin{subarray}{c} j_0+j_1+\cdots +j_{2n}=m \\ j_0,j_1,\ldots ,j_{2n} \geq 0 \end{subarray}} \zeta^{\s} ( \{ 2 \}^{j_0} , 3, \{ 2 \}^{j_1} ,1,\ldots , 3, \{ 2 \}^{j_{2n-1}} ,1 ,\{ 2 \}^{j_{2n}} ) \in \Q \cdot \pi^{4n+2m} \]
for $m=0,1$ (see \cite{H, BB, zlob, ako, muneta1} for example). Therefore the right-hand side of \eqref{c} is expressed as a rational multiple of $\pi^{4n+4}$ and we conclude the theorem.
\end{proof}
Next we prove Theorem \ref{mainthm}. For integers $a,b,c>0$ and $i,j,k\geq 0$, we put
\begin{align*}
A_{i,j} &= (z_a z_b)^i z_c (z_a z_b)^j, \\
B_{i,j} &= (z_b z_a)^i z_c (z_b z_a)^j z_b, \\
C_{i,j,k} &= (z_a z_b)^i z_c (z_a z_b)^j z_c (z_a z_b)^k, \\
D_{i,j,k} &= (z_a z_b)^i z_c (z_a z_b)^j z_a z_c (z_b z_a)^k z_b, \\
E_{i,j,k} &= (z_b z_a)^i z_c (z_b z_a)^j z_b z_c (z_a z_b)^k, \\
F_{i,j,k} &= (z_b z_a)^i z_c (z_b z_a)^j z_c (z_b z_a)^k z_b, 
\end{align*}
where $z_l=x^{l-1}y\ (l>0)$. By the definition of the product $\shc$, we obtain the following identities:  
\begin{align}\hspace{-5pt}\left.
\begin{aligned}
z_c \ \shc \ (z_a z_b)^n = & \sum_{i+j=n} A_{i,j} + \sum_{i+j=n-1} z_a B_{i,j},\\
z_c \ \shc \ z_b (z_a z_b)^n = & \sum_{i+j=n} z_b A_{i,j} +\sum_{i+j=n} B_{i,j}, \\
z_c^2 \ \shc \ (z_a z_b)^n = &\sum_{i+j+k=n} C_{i,j,k } + \sum_{i+j+k=n-1} D_{i,j,k} \\
&+\sum_{i+j+k=n-1}z_a E_{i,j,k} +\sum_{i+j+k=n-1}z_a F_{i,j,k}, \\
z_c^2 \ \shc \ z_b (z_a z_b)^n =& \sum_{i+j+k=n} z_b C_{i,j,k}+\sum_{i+j+k=n-1} z_b D_{i,j,k} \\
&+\sum_{i+j+k=n}E_{i,j,k}+\sum_{i+j+k=n}F_{i,j,k}
\end{aligned}
\right\}\label{mnm}
\end{align}
for $n\geq 0$.
\begin{proof}[Proof of Theorem \ref{mainthm}]
The proof goes by induction on $n$ such as $(\alpha_0) ,(\beta_0) \Rightarrow (\alpha_1) \Rightarrow (\beta_1) \Rightarrow (\alpha_2) \Rightarrow \cdots$. We find that the identities $(\alpha_0)$ and $(\beta_0)$ hold by simple calculation. Assuming that it has been proved up to $(\beta_{n-1})$, we prove ($\alpha_n$). The key identity is 
\begin{equation}
d(z_{k_1} \cdots z_{k_n} ) = \sum_{i=1}^n z_{k_1+\cdots +k_i} d( z_{k_{i+1}} \cdots z_{k_n} ),\label{eqn:star}
\end{equation}
where $z_{k_{i+1}} \cdots z_{k_n} =1$ if $i=n$. Using this key identity, we obtain
\begin{align*}
&\sum_{i+j+k=n} d(C_{i,j,k}) \\
&= \sum_{i+j+k=n} \sum_{h=1}^i z_{(a+b)h} d( C_{i-h,j,k})+\sum_{i+j+k=n} \sum_{h=0}^{i-1} z_{(a+b)h+a} d(z_b C_{i-h-1,j,k}) \\
&\ + \sum_{i+j+k=n} \sum_{h=0}^j z_{(a+b)(h+i)+c} d(A_{j-h,k}) + \sum_{i+j+k=n} \sum_{h=0}^{j-1} z_{(a+b)(h+i)+c+a} d(z_b A_{j-h-1,k})  \\
&\ + \sum_{i+j+k=n} \sum_{h=0}^k z_{(a+b)(h+i+j)+2c} d((z_a z_b)^{k-h}) + \sum_{i+j+k=n} \sum_{h=0}^{k-1} z_{(a+b)(h+i+j)+2c+a} d(z_b (z_a z_b)^{k-h-1}) \\
&= \sum_{h+i+j+k=n-1} z_{(a+b)(h+1)} d(C_{i,j,k})+\sum_{h+i+j+k=n-1}z_{(a+b)h+a} d(z_b C_{i,j,k}) \\
&\ +\sum_{h+i+j+k=n}  z_{(a+b)(h+i)+c} d(A_{j,k})+\sum_{h+i+j+k=n-1} z_{(a+b)(h+i)+a+c} d(z_b A_{j,k}) \\
&\ +\sum_{h+i+j+k=n} z_{(a+b)(h+i+j)+2c} d((z_az_b)^k)+\sum_{h+i+j+k=n-1} z_{(a+b)(h+i+j)+a+2c} d(z_b (z_az_b)^k) \\
&=\sum_{h+i+j+k=n-1} z_{(a+b)(h+1)} d(C_{i,j,k}) + \sum_{h+i+j+k=n-1}z_{(a+b)h+a} d(z_b C_{i,j,k}) \\
&\ +\sum_{i+j+k=n}  (i+1) z_{(a+b)i+c} d(A_{j,k})
+\sum_{i+j+k=n-1} (i+1) z_{(a+b)i+a+c} d(z_b A_{j,k}) \\
&\ +\sum_{j+k=n} \binom{j+2}{2} z_{(a+b)j+2c} d((z_az_b)^k)+\sum_{j+k=n-1} \binom{j+2}{2} z_{(a+b)j+a+2c} d(z_b (z_az_b)^k).
\end{align*}
In the same way, we find
\begin{align*}
&\sum_{i+j+k=n-1} d(D_{i,j,k}) \\
&= \sum_{h+i+j+k=n-2} z_{(a+b)(h+1)} d(D_{i,j,k}) + \sum_{h+i+j+k=n-2}z_{(a+b)h+a} d(z_bD_{i,j,k}) \\
&\ +\sum_{i+j+k=n-1} (i+1) z_{(a+b)i+c} d(z_a B_{j,k}) +\sum_{i+j+k=n-1}(i+1) z_{(a+b)i+a+c} d(B_{j,k}) \\
&\ +\sum_{j+k=n} \binom{j+1}{2} z_{(a+b)j+2c} d((z_a z_b)^k )+\sum_{j+k=n-1} \binom{j+2}{2}z_{(a+b)j+a+2c} d(z_b(z_a z_b)^k ), 
\end{align*}
\begin{align*}
& \sum_{i+j+k=n-1} d(z_a E_{i,j,k}) \\
&= \sum_{h+i+j+k=n-2} z_{(a+b)(h+1)} d(z_aE_{i,j,k}) + \sum_{h+i+j+k=n-1} z_{(a+b)h+a} d(E_{i,j,k}) \\
&\ +\sum_{i+j+k=n} i z_{(a+b)i+c} d(A_{j,k}) +\sum_{i+j+k=n-1}(i+1) z_{(a+b)i+a+c} d(z_b A_{j,k}) \\
&\ +\sum_{j+k=n} \binom{j+1}{2} z_{(a+b)j+2c} d((z_a z_b)^k) + \sum_{j+k=n-1} \binom{j+1}{2} z_{(a+b)j+a+2c} d(z_b(z_a z_b)^k) 
\end{align*}
and
\begin{align*}
&\sum_{i+j+k=n-1} d(z_a F_{i,j,k}) \\
&= \sum_{h+i+j+k=n-2}  z_{(a+b)(h+1)} d(z_a F_{i,j,k}) + \sum_{h+i+j+k=n-1}  z_{(a+b)h+a} d(F_{i,j,k}) \\
&\ +\sum_{i+j+k=n-1} i z_{(a+b)i+c} d(z_a B_{j,k}) +\sum_{i+j+k=n-1} (i+1) z_{(a+b)i+a+c} d(B_{j,k}) \\
&\ + \sum_{j+k=n} \binom{j+1}{2} z_{(a+b)j+2c} d((z_az_b)^k)+\sum_{j+k=n-1} \binom{j+2}{2}z_{(a+b)j+a+2c} d(z_b (z_az_b)^k). 
\end{align*}
These four identities add up to the left-hand side of ($\alpha_n$) because of \eqref{mnm}. 
Therefore, again using \eqref{mnm}, we obtain
\begin{align}
(\mbox{LHS of ($\alpha_n$)})
&= \sum_{j+k=n-1} z_{(a+b)(j+1)} d(z_c^2 \ \shc \ (z_a z_b )^k ) \label{eqq1} \\
&\ +\sum_{j+k=n-1}z_{(a+b)j+a} d(z_c^2 \ \shc \ z_b (z_a z_b )^k) \label{eqq2} \\
&\ +\sum_{j+k=n} (2j+1) z_{(a+b)j+c} d(z_c \ \shc \ (z_a z_b )^k ) \label{eqq3} \\
&\ +\sum_{j+k=n-1} (2j+2) z_{(a+b)j+a+c} d(z_c \ \shc \ z_b (z_a z_b )^k ) \label{eqq4} \\
&\ +\sum_{j+k=n} \binom{2j+2}{2} z_{(a+b)j+2c} d( (z_a z_b )^k ) \label{eqq5} \\
&\ +\sum_{j+k=n-1} \binom{2j+3}{2} z_{(a+b)j+a+2c} d(z_b (z_a z_b )^k )  \label{eqq6}.
\end{align}
\noindent So, it is sufficient to show that the right-hand side of the above identity equals to the right-hand side of ($\alpha_n$). 

First we have
\begin{align*}
&\sum_{j+k=n} d(A_{j,k}) \\
&=\sum_{j+k=n}\sum_{i=1}^j z_{(a+b)i}d(A_{j-i,k})+\sum_{j+k=n}\sum_{i=0}^{j-1}z_{(a+b)i+a}d(z_bA_{j-i-1,k}) \\
&\ +\sum_{j+k=n}\sum_{i=0}^k z_{(a+b)(i+j)+c}d((z_az_b)^{k-j})+\sum_{j+k=n}\sum_{i=0}^{k-1}z_{(a+b)(i+j)+a+c}d(z_b(z_az_b)^{k-i-1}) \\
&= \sum_{i+j+k=n-1} z_{(a+b)(i+1)} d(A_{j,k})+\sum_{i+j+k=n-1} z_{(a+b)i+a} d(z_b A_{j,k}) \\
&\ + \sum_{i+j+k=n} z_{(a+b)(i+j)+c} d((z_a z_b)^k)+ \sum_{i+j+k=n-1} z_{(a+b)(i+j)+a+c} d(z_b (z_a z_b)^k), \\
&= \sum_{i+j+k=n-1} z_{(a+b)(i+1)} d(A_{j,k})+\sum_{i+j+k=n-1} z_{(a+b)i+a} d(z_b A_{j,k}) \\
&\ + \sum_{i+j=n} (i+1)z_{(a+b)i+c} d((z_a z_b)^j)+ \sum_{i+j=n-1} (i+1)z_{(a+b)i+a+c} d(z_b (z_a z_b)^j).
\end{align*}
In the same way, we have
\begin{align*}
&\sum_{j+k=n-1} d(z_a B_{j,k}) \\
&=  \sum_{i+j+k=n-2} z_{(a+b)(i+1)} d(z_a B_{j,k})+\sum_{i+j+k=n-1} z_{(a+b)i+a} d(B_{j,k}) \\
&\ + \sum_{i+j=n} iz_{(a+b)i+c} d((z_a z_b)^j)+ \sum_{i+j=n-1} (i+1)z_{(a+b)i+a+c} d(z_b (z_a z_b)^j). 
\end{align*}
Using the above two identities and \eqref{mnm}, we obtain
\begin{align*} 
&d(z_c \ \shc \ (z_a z_b)^n) \\
&=\sum_{j+k=n} d(A_{j,k}) + \sum_{j+k=n-1} d(z_a B_{j,k}) \\
&= \sum_{i+j=n-1} z_{(a+b)(i+1)} d(z_c \ \shc \ (z_a z_b)^j )+\sum_{i+j=n-1} z_{(a+b)i+a} d(z_c \ \shc \ z_b (z_a z_b)^j) \\
&\ + \sum_{i+j=n} (2i+1) z_{(a+b)i+c} d((z_a z_b)^j)+ \sum_{i+j=n-1} (2i+2) z_{(a+b)i+a+c} d(z_b (z_a z_b)^j )
 \end{align*}
\noindent
for $n \geq 0$. By this identity and the harmonic product rule (\ref{har}),
we write the first term of the right-hand side of $(\alpha_n)$ (divided by the coefficient $2$) 
as
\begin{align}
& \sum_{j+k=n} d(z_c \ \shc \ (z_az_b)^j ) \ast z_{(a+b)k+c} \nonumber \\
&= \sum_{i+j+k=n-1} z_{(a+b)(i+1)} \left\{ d(z_c \ \shc \ (z_a z_b)^j) \ast z_{(a+b)k+c } \right\}  \label{eq1} \\
&\ + \sum_{j+k=n} (j+1) z_{(a+b)j+c} d(z_c \ \shc \ (z_az_b)^k)  \label{eq3} \\
&\ + \sum_{i+j+k=n-1} z_{(a+b)i+a} \left\{ d(z_c \ \shc \ z_b (z_az_b)^j )\ast z_{(a+b)k+c} \right\} \label{eq2} \\
&\ + \sum_{j+k=n-1} (j+1) z_{(a+b)j+a+c} d(z_c \ \shc \ z_b (z_az_b)^k) \label{eq5} \\
&\ + \sum_{i+j+k=n} (2i+1) z_{(a+b)i+c} \left\{ d((z_az_b)^j)\ast z_{(a+b)k+c} \right\} \label{eq4}  \\
&\ + \sum_{j+k=n} (j+1)^2 z_{(a+b)j+2c} d((z_a z_b)^k) \label{eq7} \\
&\ + \sum_{i+j+k=n-1} (2i+2) z_{(a+b)i+a+c} \left\{ d(z_b (z_az_b)^j) \ast z_{(a+b)k+c} \right\} \label{eq6}  \\
&\ + \sum_{j+k=n-1} (j+1)(j+2) z_{(a+b)j+a+2c} d(z_b (z_az_b)^k) \label{eq8} .
\end{align}
\noindent 

Note that the key identity \eqref{eqn:star} shows 
\begin{equation}
d(z_{a+b}^l)=\sum_{i+k=l-1} z_{(a+b)(i+1)} d(z_{a+b}^k) \label{eqn:a}
\end{equation}
and
\begin{equation}
d((z_a z_b)^l)=\sum_{i+j=l-1} z_{(a+b)i+a} d(z_b (z_a z_b)^j ) + \sum_{i+j=l-1} z_{(a+b)(i+1)} d((z_a z_b)^j) \label{eqn:b}
\end{equation}
for $l\geq 1$. 
By using the $\shc$-product rule \eqref{eqn:ast} and the identity \eqref{eqn:a}, the second term of the right-hand side of $(\alpha_n)$ 
is calculated as
\begin{align}
& \sum_{j+k=n} (z_c^2 \ \shc \ (z_a z_b)^j)\ast d(z_{a+b}^k) \label{second} \\
&=z_c^2\ast d(z_{a+b}^n)+\sum_{\begin{subarray}{c}j+k=n \\ j,k\geq 1 \end{subarray}}(z_c^2\ \shc \ (z_az_b)^j)\ast d(z_{a+b}^k)+z_c^2\ \shc \ (z_az_b)^n \nonumber \\
&=\sum_{j+k=n-1}z_c^2\ast z_{(a+b)(j+1)}d(z_{a+b}^k) \nonumber \\
&\ +\sum_{\begin{subarray}{c}i+j+k=n-1 \\ k\geq 1 \end{subarray}}\left\{z_c(z_c\ \shc \ (z_az_b)^k)+z_a(z_c^2\ \shc \ z_b(z_az_b)^{k-1})\right\}\ast z_{(a+b)(i+1)}d(z_{a+b}^j) \nonumber \\
&\ +z_c(z_c\ \shc \ (z_az_b)^n) \label{eqn9} \\
&\ +z_a(z_c^2\ \shc \ z_b(z_az_b)^{n-1}). \label{eqn10}
\end{align}
Expanding the first and the second terms of the right-hand side by the harmonic product rule \eqref{har}, we have
\begin{align}
& \sum_{j+k=n-1}z_c^2\ast z_{(a+b)(j+1)}d(z_{a+b}^k) \nonumber \\
&= z_c\left(z_c\ast d(z_{a+b}^n)\right) \label{eqn1} \\
&\ +\sum_{j+k=n-1}z_{(a+b)(j+1)}\left\{z_c^2\ast d(z_{a+b}^k)\right\} \label{eqn2} \\
&\ +\sum_{j+k=n-1}z_{(a+b)(j+1)+c}\left\{z_c\ast d(z_{a+b}^k)\right\} \label{eqn3}
\end{align}
and
\begin{align}
& \sum_{\begin{subarray}{c}i+j+k=n-1 \\ k\geq 1 \end{subarray}}\left\{z_c(z_c\ \shc \ (z_az_b)^k)+z_a(z_c^2\ \shc \ z_b(z_az_b)^{k-1})\right\}\ast z_{(a+b)(i+1)}d(z_{a+b}^j) \nonumber \\
&= \sum_{\begin{subarray}{c}i+j+k=n-1 \\ k\geq 1 \end{subarray}}z_{(a+b)(i+1)}\left\{(z_c^2\ \shc \ (z_az_b)^k)\ast d(z_{a+b}^j)\right\} \label{eqn4} \\
&\ +\sum_{\begin{subarray}{c}j+k=n \\ j,k\geq 1 \end{subarray}}z_c\left\{(z_c\ \shc \ (z_az_b)^j) \ast d(z_{a+b}^k)\right\} \label{eqn5} \\
&\ +\sum_{\begin{subarray}{c}i+j+k=n-1 \\ k\geq 1 \end{subarray}}z_{(a+b)(i+1)+c}\left\{(z_c\ \shc \ (z_az_b)^k) \ast d(z_{a+b}^j)\right\} \label{eqn6} \\
&\ +\sum_{\begin{subarray}{c}j+k=n \\ j,k\geq 1 \end{subarray}}z_a\left\{(z_c^2\ \shc \ z_b(z_az_b)^{j-1}) \ast d(z_{a+b}^k)\right\} \label{eqn7} \\
&\ +\sum_{\begin{subarray}{c}i+j+k=n-1 \\ k\geq 1 \end{subarray}}z_{(a+b)(i+1)+a}\left\{(z_c^2\ \shc \ z_b(z_az_b)^{k-1}) \ast d(z_{a+b}^j)\right\}. \label{eqnstar} 
\end{align}
We see that identities
\begin{align*}
&\eqref{eqn2}+\eqref{eqn4}=\sum_{i+j+k=n-1} z_{(a+b)(i+1)} \left\{ (z_c^2 \ \shc \ (z_az_b)^j)\ast d(z_{a+b}^k) \right\}, \\
&\eqref{eqn9}+\eqref{eqn1}+\eqref{eqn5}=\sum_{j+k=n} z_{c} \left\{ (z_c \ \shc \ (z_a z_b)^j)\ast d(z_{a+b}^k) \right\}, \\
&\eqref{eqn10}+\eqref{eqn7}=\sum_{j+k=n-1} z_{a} \left\{ (z_c^2 \ \shc \ z_b (z_a z_b)^j)\ast d(z_{a+b}^k) \right\}, \\
&\eqref{eqn3}+\eqref{eqn6}=\sum_{i+j+k=n-1} z_{(a+b)(i+1)+c} \left\{ (z_c \ \shc \ (z_a z_b)^j)\ast d(z_{a+b}^k) \right\}, \\
&\eqref{eqnstar}=\sum_{\begin{subarray}{c}i+j+k=n-1 \\ i\geq 1 \end{subarray}}z_{(a+b)i+a}\left\{(z_c^2\ \shc \ z_b(z_az_b)^j)\ast d(z_{a+b}^k)\right\}
\end{align*}
hold. Therefore we have
\begin{align}
\eqref{second}&= \sum_{i+j+k=n-1} z_{(a+b)(i+1)} \left\{ (z_c^2 \ \shc \ (z_az_b)^j)\ast d(z_{a+b}^k) \right\} \label{eq9} \\
&\ + \sum_{i+j+k=n} z_{(a+b)i+c} \left\{ (z_c \ \shc \ (z_a z_b)^j)\ast d(z_{a+b}^k) \right\} \label{eq10} \\
&\ + \sum_{i+j+k=n-1} z_{(a+b)i+a} \left\{ (z_c^2 \ \shc \ z_b (z_a z_b)^j)\ast d(z_{a+b}^k) \right\} \label{eq11}.
\end{align}

By \eqref{eqn:b} and the harmonic product rule \eqref{har}, the third term of the right-hand side of $(\alpha_n)$ (divided by the coefficient $-4$) is calculated as
\begin{align}
& \sum_{i+j+k=n} d((z_a z_b)^i)\ast z_{(a+b)j+c} z_{(a+b)k+c} \label{third} \\
&= \sum_{j+k=n}z_{(a+b)j+c}z_{(a+b)k+c}+\sum_{\begin{subarray}{c}i+j+k=n \\ i\geq 1 \end{subarray}}d((z_az_b)^i)\ast z_{(a+b)j+c}z_{(a+b)k+c} \nonumber \\
&= \sum_{j+k=n}z_{(a+b)j+c}z_{(a+b)k+c} \nonumber \\
&\ +\sum_{h+i+j+k=n-1}z_{(a+b)h+a}d(z_b(z_az_b)^i)\ast z_{(a+b)j+c}z_{(a+b)k+c} \nonumber \\
&\ +\sum_{h+i+j+k=n-1}z_{(a+b)(h+1)}d((z_az_b)^i)\ast z_{(a+b)j+c}z_{(a+b)k+c} \nonumber \\
&= \sum_{j+k=n}z_{(a+b)j+c}z_{(a+b)k+c} \label{eqnast} \\
&\ +\sum_{h+i+j+k=n-1}\Bigl\{z_{(a+b)h+a}\left(d(z_b(z_az_b)^i)\ast z_{(a+b)j+c}z_{(a+b)k+c} \right) \nonumber \\
&\qquad\qquad\qquad\quad +z_{(a+b)j+c}\left(z_{(a+b)h+c}d(z_b(z_az_b)^i)\ast z_{(a+b)k+c} \right) \label{eqn11} \\
&\qquad\qquad\qquad\quad +z_{(a+b)(h+j)+a+c}\left(d(z_b(z_az_b)^i)\ast z_{(a+b)k+c} \right) \label{eqn12} \\
&\qquad\qquad\qquad\quad +z_{(a+b)(h+1)}\left(d((z_az_b)^i)\ast z_{(a+b)j+c}z_{(a+b)k+c} \right) \nonumber \\
&\qquad\qquad\qquad\quad +z_{(a+b)j+c}\left(z_{(a+b)(h+1)}d((z_az_b)^i)\ast z_{(a+b)k+c} \right) \label{eqn13} \\
&\qquad\qquad\qquad\quad +z_{(a+b)(h+j+1)+c}\left(d((z_az_b)^i)\ast z_{(a+b)k+c} \right)\Bigr\}. \label{eqn14}
\end{align}
By \eqref{eqn:b}, 
\[ \eqref{eqn11}+\eqref{eqn13}=\sum_{\begin{subarray}{c}i+j+k=n \\ j\geq 1 \end{subarray}}z_{(a+b)i+c}\left(d((z_az_b)^j)\ast z_{(a+b)k+c}\right). \]
Also we have
\[ \eqref{eqn14}=\sum_{i+j+k=n}iz_{(a+b)i+c}\left(d((z_az_b)^j)\ast z_{(a+b)k+c}\right). \]
Hence we find
\[ \eqref{eqnast}+\eqref{eqn11}+\eqref{eqn13}+\eqref{eqn14}= \sum_{i+j+k=n} (i+1) z_{(a+b)i+c} \left\{ d((z_az_b)^j)\ast z_{(a+b)k+c} \right\}. \]
Since
\[ \eqref{eqn12}= \sum_{i+j+k=n-1} (i+1) z_{(a+b)i+a+c} \left\{ d(z_b (z_a z_b)^j) \ast z_{(a+b)k+c}\right\}, \]
we have
\begin{align}
\eqref{third}&= \sum_{h+i+j+k=n-1} z_{(a+b)(h+1)} \left\{ d((z_a z_b)^i) \ast z_{(a+b)j+c} z_{(a+b)k+c} \right\} \label{eq12} \\
&\ + \sum_{h+i+j+k=n-1} z_{(a+b)h+a} \left\{ d(z_b (z_a z_b)^i)\ast z_{(a+b)j+c} z_{(a+b)k+c} \right\} \label{eq13} \\
&\ + \sum_{i+j+k=n} (i+1) z_{(a+b)i+c} \left\{ d((z_az_b)^j)\ast z_{(a+b)k+c} \right\} \label{eq14} \\
&\ + \sum_{i+j+k=n-1} (i+1) z_{(a+b)i+a+c} \left\{ d(z_b (z_a z_b)^j) \ast  z_{(a+b)k+c} \right\}. \label{eq15} 
\end{align}

By \eqref{eqn:b} and the harmonic product rule \eqref{har}, the fourth term of the right-hand side of $(\alpha_n)$ (divided by the coefficient $-1$) is calculated as
\begin{align}
& \sum_{j+k=n} d((z_a z_b )^j)\ast z_{(a+b)k+2c} \label{fourth} \\
&= z_{(a+b)n+2c}+\sum_{\begin{subarray}{c}j+k=n \\ j\geq 1 \end{subarray}}d((z_az_b)^j)\ast z_{(a+b)k+2c} \nonumber \\
&= z_{(a+b)n+2c}+\sum_{i+j+k=n-1}z_{(a+b)i+a}d(z_b(z_az_b)^j)\ast z_{(a+b)k+2c} \nonumber \\
&\ +\sum_{i+j+k=n-1}z_{(a+b)(i+1)}d((z_az_b)^j)\ast z_{(a+b)k+2c} \nonumber \\
&= z_{(a+b)n+2c} \label{eqn15} \\
&\ +\sum_{i+j+k=n-1}\Bigl\{z_{(a+b)i+a}\left(d(z_b(z_az_b)^j)\ast z_{(a+b)k+2c}\right) \nonumber \\
&\qquad\qquad\qquad +z_{(a+b)k+2c}z_{(a+b)i+a}d(z_b(z_az_b)^j) \label{eqn16} \\
&\qquad\qquad\qquad +z_{(a+b)(i+k)+a+2c}d(z_b(z_az_b)^j) \label{eqn17} \\
&\qquad\qquad\qquad +z_{(a+b)(i+1)}\left(d((z_az_b)^j)\ast z_{(a+b)k+2c}\right) \nonumber \\
&\qquad\qquad\qquad +z_{(a+b)k+2c}z_{(a+b)(i+1)}d((z_az_b)^j) \label{eqn18} \\
&\qquad\qquad\qquad +z_{(a+b)(i+k+1)+2c}d((z_az_b)^j)\Bigr\}.\label{eqn19}
\end{align}
By \eqref{eqn:b}, 
\[ \eqref{eqn16}+\eqref{eqn18}=\sum_{\begin{subarray}{c}j+k=n \\ j\geq 1 \end{subarray}}z_{(a+b)k+2c}d((z_az_b)^j). \]
Also we have
\[ \eqref{eqn19}=\sum_{j+k=n}kz_{(a+b)k+2c}d((z_az_b)^j). \]
Hence we find
\[ \eqref{eqn15}+\eqref{eqn16}+\eqref{eqn18}+\eqref{eqn19}=\sum_{j+k=n}  (j+1) z_{(a+b)j+2c} d((z_a z_b)^k). \]
Since
\[ \eqref{eqn17}= \sum_{j+k=n-1}  (j+1) z_{(a+b)j+a+2c} d(z_b(z_az_b)^k), \]
we have
\begin{align}
\eqref{fourth}& = \sum_{i+j+k=n-1} z_{(a+b)(i+1)} \left\{  d((z_a z_b)^j)\ast z_{(a+b)k+2c} \right\} \label{eq16} \\
&\ + \sum_{i+j+k=n-1} z_{(a+b)i+a} \left\{  d(z_b (z_a z_b)^j) \ast z_{(a+b)k+2c} \right\} \label{eq17} \\
&\ + \sum_{j+k=n}  (j+1) z_{(a+b)j+2c} d((z_a z_b)^k) \label{eq18} \\
&\ + \sum_{j+k=n-1}  (j+1) z_{(a+b)j+a+2c} d(z_b(z_az_b)^k). \label{eq19}
\end{align}

Following identity is shown by \cite[Theorem 7 (5)]{muneta1} (where the map $d$ is denoted by $S$) and \eqref{mnm}: 
\begin{equation}
d(z_c\ \shc \ (z_az_b)^n)=2\sum_{j+k=n}d((z_az_b)^j)\ast z_{(a+b)k+c}-\sum_{j+k=n}\left(z_c\ \shc \ (z_az_b)^j\right)\ast d(z_{a+b}^k) \label{mun}
\end{equation}
for $n\geq 0$. 
By induction hypothesis and 
\eqref{mun}, we have
\begin{align*}
(\ref{eqq1})&=2\times (\ref{eq1})+(\ref{eq9})-4 \times (\ref{eq12}) -(\ref{eq16}) \ \ \ (\mbox{by} \ (\alpha_0),\ldots ,(\alpha_{n-1})), \\
(\ref{eqq2})&=2\times (\ref{eq2})+(\ref{eq11})-4 \times (\ref{eq13}) -(\ref{eq17}) \ \ \ (\mbox{by} \ (\beta_0),\ldots ,(\beta_{n-1})), \\
(\ref{eqq3})&=2\times (\ref{eq3})+2\times (\ref{eq4}) +(\ref{eq10})-4\times (\ref{eq14}) \ \ \ (\mbox{by}\  \eqref{mun}). 
\end{align*}
Also we immediately find that
\begin{align*}
(\ref{eqq4})&=2\times (\ref{eq5}), \\
(\ref{eqq5})&=2\times (\ref{eq7})-(\ref{eq18}), \\
(\ref{eqq6})&=2\times (\ref{eq8})-(\ref{eq19}), \\
0 &= 2\times (\ref{eq6}) -4 \times (\ref{eq15}).
\end{align*}
We have already observed that
the left(resp. right)-hand sides of these seven equations add up to the left(resp. right)-hand side of the identity $(\alpha_n)$. Therefore we prove $(\alpha_n)$.
In the same way, $(\beta_n)$ can be proved by using the induction hypothesis and $(\alpha_n)$ for $n$.
This completes the proof of Theorem \ref{mainthm}. 
\end{proof}

\section{Conjectures}\label{sec4}
We discover more general formulas when we search relations among MZV's and MZSV's based on numerical experiments. For non-negative integers $j_0, \ldots ,j_{2n}$, put
\[ \overline\Z (j_0,j_1, \ldots ,j_{2n} ):=\zeta^{\s} (\{ 2 \}^{j_0},3,\{ 2 \}^{j_1} ,1,\{ 2 \}^{j_2}, \ldots ,3,\{ 2\}^{j_{2n-1}},1,\{ 2\}^{j_{2n}}). \]
We write two operators as
\[ (j_0,j_1, \ldots ,j_n )_{+} = (j_0,j_1, \ldots ,j_{n-1} ,j_n ,0 ), \]
\[ (j_0,j_1, \ldots ,j_n )^{+} = (j_0,j_1, \ldots ,j_{n-1} ,j_{n}+1 ). \]
Let us denote by $\sym_n$ the symmetric group of degree $n$. For a vector $S=(j_0 , \ldots , j_n )$ and $\sigma\in\sym_{n+1}$, we define 
\[ \sigma S= (j_{\sigma (0) } , \ldots ,j_{\sigma (n)} ). \]
\begin{conj} \label{conj1} Let $n$ be a positive integer, $j_0 , \ldots , j_{2n-1}$ non-negative integers. Put $m=j_0 + \cdots +j_{2n-1}$. For a vector $S=(j_0, \ldots ,j_{2n-1})$, we have
\begin{align*} \sum_{ \sigma \in \sym_{2n} } \overline{\Z}((\sigma S)_{+}) \stackrel{?}{\in} \Q \cdot  \pi^{2m+4n}.
\end{align*}
\end{conj}
\begin{exmp}\label{ex1}
For the vector $S=(0,0)$, Conjecture \ref{conj1} reads
\[ \zeta^{\s}(3,1)\in\Q\cdot\pi^4, \]
which is shown in many ways. For the vector $S=(1,0)$, Conjecture \ref{conj1} reads
\[ \zeta^{\s}(2,3,1)+\zeta^{\s}(3,2,1)\stackrel{?}{\in}\Q\cdot\pi^6, \]
which is verified in section \ref{sec5}. For the vector $S=(1,0,0,0)$, Conjecture \ref{conj1} reads
\[ \zeta^{\s}(2,3,1,3,1)+\zeta^{\s}(3,2,1,3,1)+\zeta^{\s}(3,1,2,3,1)+\zeta^{\s}(3,1,3,2,1)\stackrel{?}{\in}\Q\cdot\pi^{10}. \]
\end{exmp}
\begin{conj} \label{conj11} Let $n,j_0,\ldots , j_{2n}$ be non-negative integers. Put $m=j_0 + \cdots +j_{2n} $.  For a vector $S=(j_0,\ldots , j_{2n})$, we have
\begin{align*} 
\sum_{\sigma \in \sym_{2n+1}} \overline{\Z}( (\sigma S)^{+}) \stackrel{?}{\in} \Q \cdot  \pi^{2m+4n+2}.
\end{align*} 
\end{conj}
\begin{exmp}\label{ex2}
For the vector $S=(0)$, Conjecture \ref{conj11} reads
\[ \zeta^{\s}(2)\in\Q\cdot\pi^2, \]
which is well known. For the vector $S=(0,0,0)$, Conjecture \ref{conj11} reads
\[ \zeta^{\s}(3,1,2)\stackrel{?}{\in}\Q\cdot\pi^6, \]
which is verified in section \ref{sec5}. For the vector $S=(1,0,0)$, Conjecture \ref{conj11} reads
\[ \zeta^{\s}(2,3,1,2)+\zeta^{\s}(3,2,1,2)+\zeta^{\s}(3,1,2,2)\stackrel{?}{\in}\Q\cdot\pi^8. \]
\end{exmp}
\noindent We expect that these conjectures are primitive, i.e., partial sums of the left-hand side of our conjectures cannot be rational multiples of powers of $\pi^2$. Conjecture \ref{conj1} reduces to $\zeta^{\s} ( \{ 3,1 \}^n) \in \Q \cdot \pi^{4n} $ if $m=0$, and Conjecture \ref{conj11} to $\zeta^{\s} (\{ 2 \}^m) \in \Q \cdot \pi^{2m} $ if $n=0$. Both have been already proved (for example, see \cite{zlob,muneta1}).

Moreover, we experimentally verify the following formulas. The formulas can be proved algebraically for lower weights by using known relations such as extended double shuffle relation (\cite{IKZ}). Some are discussed in the next section. But in general, we have not proved them yet. 
\begin{conj}\label{tasaka} A) For any integers $n,m \geq 0$, we have
\begin{align*}  \zeta^{\s} ( \{ 2 \}^n ,3,\{ 2 \}^m ,1 ) +\zeta^{\s} ( \{ 2 \}^m ,3,\{ 2 \}^n ,1 ) \stackrel{?}{=} \zeta^{\s} ( \{ 2 \}^{n+1} ) \zeta^{\s} ( \{ 2 \}^{m+1} ). 
\end{align*}
B) For any integer $n \geq 0$, we have
\begin{align*} (2n+1) \zeta^{\s} ( \{ 3,1 \}^n ,2 ) \stackrel{?}{=} \sum_{j+k=n} \zeta^{\s} (\{ 3,1 \}^j ) \zeta^{\s} (\{ 2 \}^{2k+1} ). 
\end{align*} 
C) For any integer $n\geq 1$, we have 
\begin{align*} \sum_{\begin{subarray}{c} j_0+j_1+\cdots +j_{2n-1}=1 \\ j_0,j_1,\ldots ,j_{2n-1} \geq 0 \end{subarray}} \zeta^{\s} ( \{ 2 \}^{j_0} , 3, \{ 2 \}^{j_1} ,1,\ldots , 3, \{ 2 \}^{j_{2n-1}} ,1 ) \\
\stackrel{?}{=}  \sum_{j+k=n-1} \zeta^{\s} ( \{ 3,1 \}^j ,2) \zeta^{\s} (\{ 2 \}^{2k+2}). 
\end{align*}
\end{conj}

\section{Remarks}\label{sec5}
Instead of the assertion A) of Conjecture \ref{tasaka}, we prove the following weaker identity:
\begin{prop} For any non-negative integer $n$, we have
\[ 2\sum_{i+j=n} \zeta^{\s} ( \{ 2 \}^i , 3, \{ 2 \}^j ,1 ) = \sum_{i+j=n} \zeta^{\s} (\{ 2 \}^{i+1} ) \zeta^{\s} (\{ 2 \}^{j+1} ). \] 
\end{prop}

\begin{proof}
Using the cyclic sum formula (\cite{ohno}) for the index $(3,\underbrace{2,\ldots ,2}_{n-2})$, we have
\[ \sum_{i+j=n-2} \zeta^{\s} ( \{ 2 \}^i  ,3, \{ 2 \}^j ,1 ) + \zeta^{\s} ( \{ 2 \}^n ) = (2n-1) \zeta (2n) \]
for $n\geq 1$. We also find 
\[ \zeta^{\s} ( \{ 2 \}^n )= 2 (1-2^{1-2n}) \zeta (2n) \]
(see \cite{aoki,zlob}). Because of the Euler's formula $ \zeta (2n)= (-1)^{n+1} \frac{B_{2n}(2\pi)^{2n}}{2 (2n)!} $ (where $B_j$ denotes the $j$-th Bernoulli number), it is enough to show
\[ (1-2n) \frac{ B_{2n}}{(2n)!} =\sum_{i+j=n} (1-2^{1-2i}) (1-2^{1-2j}) \frac{ B_{2i}B_{2j}}{(2i)!(2j)!}. \]
This equality holds since the generating function of the left-hand side coincides with that of the right-hand side, that is, 
\begin{align*} \sum_{n=0}^{\infty} (1-n) \frac{ B_{n}}{n!} t^n &= \frac{t}{e^t-1} -t\frac{d}{dt} \left( \frac{t}{e^t-1} \right) \\
&= \left( \frac{t}{e^t-1} - 2\frac{t/2}{e^{t/2}-1} \right)^2 \\
&= \left( \sum_{n=0}^{\infty} (1-2^{1-n}) \frac{ B_{n}}{n!} t^n \right)^2. 
\end{align*}
\end{proof} 

Lastly, we give a verification of A) and B) of Conjecture \ref{tasaka} in the first non-trivial case (in the case of $m=n=1$) by using the extended double shuffle relation proved in \cite{IKZ}. 

The extended double shuffle relation states as follows (see \cite{IKZ} for details).
Let $\sh:\h \times \h \rightarrow \h$ be the $\Q$-bilinear map defined, for any words $w,w' \in \h$ and $u,v \in \{ x,y \}$, by
\[ 1 \ \sh \ w=w \ \sh \ 1=w \]
\noindent and
\[ uw \ \sh \ vw' =u(w \ \sh \ vw')+v(uw \ \sh \ w'). \]
\noindent The product $\sh$ is commutative and associative on $\h$, which is called the shuffle product. We immediately find that $\h^1 \ \sh \ \h^1 \subset \h^1$ and $\h^0 \ \sh \ \h^0 \subset \h^0$ hold. We denote by $\h_{\sh}^1$ (resp. $\h_{\sh}^0$) the commutative algebra $\h^1$ (resp. $\h^0$) equipped with the shuffle product. It is known that
\[ \h_{\sh}^1 \cong \h_{\sh}^0 [y] \]
\noindent holds (see \cite[Theorem 6.1]{reu}), which means that any $w \in \h_{\sh}^1 $ can be expressed as
\[ w= w_0 + w_1\ \sh \ y + w_2\ \sh \ y^{\sh 2} + \cdots + w_n \ \sh \ y^{\sh n}  \quad (w_0,\ldots ,w_n \in \h^0), \]
where $y^{\sh i} =\underbrace{y \ \sh \ \cdots \ \sh \ y}_{i}$. 
We denote its constant term $w_0$ by $\reg (w)$. Then the extended double shuffle relation states as follows.
\begin{thm}[Ihara-Kaneko-Zagier] For any $w_1 \in \h^1$ and any $w_0 \in \h^0$, we have
\begin{equation*} Z (\reg (w_1 \ \sh \ w_0 -w_1 \ast w_0)) = 0. \end{equation*} 
\end{thm}
\noindent By using Risa/asir, an open source general computer algebra system, we found the following two identities (of weight $6$):
\begin{align*}
& d(z_2 z_3 z_1) +d(z_3 z_2 z_1) - d (z_2 z_2) \ast d(z_2) \\
&=  -\reg (xy \ast yx^2y - xy \ \sh \ yx^2y) - \reg (x^2y \ast x^2y -x^2y \ \sh \ x^2y) \\
&- \reg (x^2y \ast xy^2 - x^2y \ \sh \ xy^2 ) +4 \reg (x^3 y \ast y^2 - x^3 y \ \sh \ y^2 ) \\
&+ 2 \reg (x^4 y\ast y-x^4y \ \sh \ y) -5 \reg (x^3y^2\ast y -x^3y^2 \ \sh \ y) \\
& +2 \reg (x^2yxy\ast y -x^2yxy \ \sh \ y) +5 \reg (xy x^2y\ast y -xyx^2y \ \sh \ y)
\end{align*} 
\noindent and
\begin{align*}
& 3d(z_3 z_1 z_2) -d(z_3 z_1)\ast d(z_2) - d(z_2 z_2 z_2) \\
&= - \reg (x^2y\ast xy^2-x^2y \ \sh \ xy^2 ) + \reg (x^2y^2\ast xy-x^2y^2 \ \sh \ xy ) \\
& +2 \reg (x^3y\ast y^2-x^3y \ \sh \ y^2 ) +\reg (x^4y\ast y-x^4y \ \sh \ y ) \\
& +2 \reg (x^3y^2\ast y-x^3y^2 \ \sh \ y ) +2 \reg (x^2yxy\ast y-x^2yxy \ \sh \ y ) \\
& +2 \reg (xyx^2y \ast y - xyx^2y \ \sh \ y).
\end{align*}
\noindent 
Since the weight is relatively small, these formulas are verified by calculating perseveringly. Applying the evaluation map $Z$ to these identities, we conclude A) and B) of Conjecture \ref{tasaka} when $m=n=1$, which also give a verification of the second formulas stated in Example \ref{ex1} and \ref{ex2}.

%
%


\end{document}